\documentclass{amsart}
\usepackage[cp1251]{inputenc}
\usepackage[english]{babel}
\usepackage{amsmath}
\usepackage{amssymb}
\usepackage{amsfonts}
\newtheorem{cor}{Corollary}
\newtheorem{thm}{Theorem}

\newtheorem{defin}{Definition}
\begin{document}

\thispagestyle{empty}

 \title[ Multilayer fields by the transformation operators method ]{  Multi-component mathematical models in multilayer fields by the  transformation operators method}%

\author{O. Yaremko, L.Simutina}

\address{Oleg Yaremko, Lidia Simutina
\newline\hphantom{iii}Penza State University,
\newline\hphantom{iii}str. Krasnaya, 40, 
\newline\hphantom{iii} 440038, Penza, Russia}
\email{yaremki@mail.ru}

\maketitle {\small

\begin{quote}
\noindent{\bf Abstract. }  The vector transform operators are investigated; these operators are used at the solution of boundary
value problems in piecewise homogeneous spherically symmetric areas. In particular,
examples of transformation operators for vector boundary value problems are given
for third vector boundary value problem in the unit circle and for the Dirichlet problem in the unit circle.

\medskip

\noindent{\bf Keywords:}
Fourier matrix transforms, integral transform, heat conduction equation,  wave equation, Poisson equation.
\end{quote} }

\emph{{Mathematics Subject Classification 2010}:{35N30 Over determined initial-boundary value problems; 35Cxx Representations of solutions; 65R10 Integral transforms, 735Qxx		Equations of mathematical physics and other areas of application }.}\\

{Penza state university, PO box 440026, Penza, Krasnaya street, 40, Russia}

\section{Introduction}

  The important mathematical models describing physical fields in multilayered piecewise homogeneous media, lead to an initial boundary value problems for partial differential equations. The coefficients of the equations are continuous for homogeneous models and the coefficients of the system equation are piecewise constant for piecewise-homogeneous models. Transformation operators allow to interpret piecewise-homogeneous physical fields as a perturbing homogeneous fields.  The authors to consider the transformation operator as a deformation of the mathematical homogeneous model into a piecewise-homogeneous model. In this work vector transformation operators are constructed, studied and used for problems' solving on the composite real line.
  
  Consider  vector Robin and Dirichlet boundary value problems over the real semiaxis.First:
\begin{equation}\left\{ \begin{matrix}
   {{u}_{yy}}+{{a}^{2}}{{u}_{xx}}=0, & x>0,y\in R,  \\
   hu(0,y)+{{{{u}'}}_{x}}(0,y)=f(y), & x>0,y\in R,  \\
\end{matrix} \right.	\end{equation}

$$f(y)={{\left( {{f}_{1}}(y)...{{f}_{n}}(y) \right)}^{T}}.
$$
Second:
\begin{equation}\left\{ \begin{matrix}
   {{{\tilde{u}}}_{i,yy}}+{{{\tilde{u}}}_{i,xx}}=0, & x>0,y\in R,  \\
   {{{\tilde{u}}}_{i}}(0,y)={{f}_{i}}(y), & x>0,y\in R,i=\overline{1,n}.  \\
\end{matrix} \right.\end{equation}

\begin{thm}If the function  $\tilde{u}(x,y)$ is a solution of the Dirichlet value problem (2) spectrum $\sigma (a)$ in right semiplane and spectrum $\sigma (h)$ in left semiplane, then the vector function
$$u(x,y)=\sum\limits_{k=1}^{n}{\int_{0}^{\infty }{{{{\tilde{u}}}_{k}}}}({{a}^{-1}}x+\varepsilon I){{e}^{ah\varepsilon }}a{{e}_{k}}d\varepsilon ,{{e}_{k}}={{\left( {{0...1}_{k}}...0 \right)}^{T}}$$
is a solution of the Robin value problem (1), were \[{{e}^{A}}=\sum\limits_{k=0}^{\infty }{{{A}^{k}}/k}!,\]I-identity matrix.
\end{thm}
	Consider  vector Dirichlet boundary value problems over the real semiaxis.First:
\begin{equation}\left\{ \begin{matrix}
   {{u}_{1,yy}}+a_{1}^{2}{{u}_{1,xx}}=0,0<x<l,y\in R,  \\
   {{u}_{2,yy}}+a_{2}^{2}{{u}_{2,xx}}=0,l<x,y\in R,  \\
   {{u}_{1}}(0,y)=f(y),y\in R,  \\
   {{u}_{1}}(l,y)={{u}_{2}}(l,y),y\in R,  \\
   {{\lambda }_{1}}{{u}_{1}}(l,y)={{\lambda }_{2}}{{u}_{2}}(l,y),y\in R.  \\
\end{matrix} \right.\end{equation}
Second:
\begin{equation}\left\{ \begin{matrix}
   {{{\tilde{u}}}_{i,yy}}+{{{\tilde{u}}}_{i,xx}}=0, & x>0,y\in R,  \\
   {{{\tilde{u}}}_{i}}(0,y)={{f}_{i}}(y), & x>0,y\in R,i=\overline{1,n}.  \\
\end{matrix} \right.\end{equation}
 Let a be a square $N\times N$matrix with N linearly independent eigenvectors, ${{q}_{i}}\,\,(i=1,\ldots ,N),$then a can be factorized as $$a=Q\Lambda {{Q}^{-1}},$$ where $Q$ is the square $N\times N$ matrix whose ${{i}^{th}}$ column is the eigenvector ${{q}_{i}}$ of $a$ and $\Lambda $ is the diagonal matrix whose diagonal elements are the corresponding eigenvalues, i.e., $${{\Lambda }_{ii}}={{\lambda }_{i}}.$$
\begin{defin}Vector shift operator is defined $${{T}_{a}}f(x)=\sum\limits_{i=1}^{n}{f}(x+\lambda _{i}^{-1}l){{Q}^{-1}}{{I}_{i}}Q,$$ where the matrix ${{I}_{i}}$ of size $N\times N$ is matrix in which all the elements are equal to 0 but one element ${{I}_{ii}}=1$;
vector contraction operator by function  : $${{U}_{\chi }}f(x)=f(x)\chi ;\quad \chi ={{a}_{1}}\lambda _{1}^{-1}{{\lambda }_{2}}a_{2}^{-1};$$
reflection operator  :$Sf(x)=f(2l-x);$
vector contraction operator by argument x with the shift -l :$${{R}_{a}}f(x)=f({{a}^{-1}}(x-l)).$$ \end{defin}
\begin{thm}If the function  $\tilde{u}(x,y)$ is a solution of the vector Dirichlet value problem (4) spectrum $\sigma ({{a}_{1}}),\sigma ({{a}_{2}})$ in right semiplane ,then the vector function
$${{u}_{1}}=\sum\limits_{k=1}^{n}{\sum\limits_{j=0}^{\infty }{(}}{{R}_{{{a}_{1}}}}-S{{R}_{{{a}_{1}}}}{{U}_{\chi }}){{T}_{{{a}_{1}}}}{{({{T}_{{{a}_{1}}}}{{U}_{\chi }}{{T}_{{{a}_{1}}}})}^{j}}{{\tilde{u}}_{k}}(x){{e}_{k}},0<x<l;$$
$${{u}_{2}}=\sum\limits_{k=1}^{n}{\sum\limits_{j=0}^{\infty }{(}}{{R}_{{{a}_{2}}}}-{{R}_{{{a}_{2}}}}{{U}_{\chi }}){{T}_{{{a}_{1}}}}{{({{T}_{{{a}_{1}}}}{{U}_{\chi }}{{T}_{{{a}_{1}}}})}^{j}}{{\tilde{u}}_{k}}(x){{e}_{k}},l<x<\infty ;$$
is a solution of the vector Dirichlet value problem (3).
\end{thm}
\begin{cor}Zeroth-order approximation solution of the problem (3)  	is given by
$u_{1}^{(0)}=\sum\limits_{k=1}^{n}{{{R}_{{{a}_{1}}}}}{{T}_{{{a}_{1}}}}{{\tilde{u}}_{k}}(x){{e}_{k}},0<x<l;$$u_{2}^{(0)}=\sum\limits_{k=1}^{n}{{{R}_{{{a}_{2}}}}}{{T}_{{{a}_{1}}}}{{\tilde{u}}_{k}}(x){{e}_{k}},l<x<\infty .$
The first order approximation 	solution of the problem (3) is of the form
$$u_{1}^{(1)}=u_{1}^{(0)}+\sum\limits_{k=1}^{n}{(}{{R}_{{{a}_{1}}}}{{T}_{{{a}_{1}}}}{{T}_{{{a}_{1}}}}-S{{R}_{{{a}_{1}}}})U\chi {{T}_{{{a}_{1}}}}{{\tilde{u}}_{k}}(x){{e}_{k}},0<x<l;$$
$$u_{2}^{(1)}=u_{2}^{(0)}+\sum\limits_{k=1}^{n}{(}{{R}_{{{a}_{2}}}}{{T}_{{{a}_{1}}}}{{T}_{{{a}_{1}}}}-{{R}_{{{a}_{2}}}})U\chi {{T}_{{{a}_{1}}}}{{\tilde{u}}_{k}}(x){{e}_{k}},l<x<\infty .$$
\end{cor}


\begin{thebibliography}{99}
\bibitem {1} O.E. Yaremko, {\it Matrix integral Fourier transforms for problems with discontinuous coefficients and transformation operators}, Reports Of Academy Of Sciences, Volume. 417, Issue 3, 2007, p. 323-325.
\bibitem {2}  Bavrin I.I., Matrosov V.L., Jaremko O. E. {\it Operators of transformationin the analysis, mathematical physics and Pattern recognition}. Moscow,
Prometheus, (2006),p 292.
\bibitem {3}  O. Yaremko, V. Selutin, N. Yaremko,  {\it The Fourier Transform with Piecewise Trigonometric Kernels and its Applications}, WSEAS transactions on mathematics, Volume 13, 2014, pp. 615-625.
\bibitem {4} O.E. Yaremko,  {\it Transformation operator and boundary value problems}, Differential Equation. Vol.40, No. 8, 2004, pp.1149-1160
\end{thebibliography}
\end{document}